\documentclass[11pt]{amsart}
\usepackage{mathrsfs}
\usepackage{amssymb,latexsym}
\usepackage{pstcol,pstricks,color}

 \numberwithin{equation}{section}

\newtheorem{theorem}{Theorem}[section]

\newtheorem{conjecture}[theorem]{Conjecture}
\newtheorem{remark}[theorem]{Remark}

\title{ Another simple proof of an identity conjectured by Lacasse }

\begin{document}
\maketitle
\begin{center}
Yidong Sun$^\dag$

%%%%\thanks{Corresponding author: Yidong Sun, sydmath@yahoo.com.cn.}%%

Department of Mathematics, Dalian Maritime University, 116026 Dalian, P.R. China\\[5pt]

{\it  Emails: $^\dag$sydmath@yahoo.com.cn }

\end{center}\vskip0.2cm

\subsection*{Abstract}
In this note, using the derangement polynomials and their umbral
representation, we give another simple proof of an identity
conjectured by Lacasse in the study of the PAC-Bayesian machine
learning theory.

\medskip

{\bf Keywords}: Derangement polynomial; Umbral operator.

\noindent {\sc 2000 Mathematics Subject Classification}: 05A15,
05A19.

{\bf \section{Introduction} }

In his thesis, Lacasse introduced the functions $\xi(n)$ and $\xi_2(n)$
in the study of the PAC-Bayesian machine learning theory, where
\begin{eqnarray*}
\xi(n)  \hskip-.22cm &=&\hskip-.22cm \sum_{k=0}^{n}\binom{n}{k}\left(\frac{k}{n}\right)^k\left(1-\frac{k}{n}\right)^{n-k}, \\
\xi_2(n)\hskip-.22cm &=&\hskip-.22cm \sum_{k=0}^{n}\sum_{j=0}^{n-k}\binom{n}{k}\binom{n-k}{j}\left(\frac{k}{n}\right)^k\left(\frac{j}{n}\right)^j\left(1-\frac{k}{n}-\frac{j}{n}\right)^{n-k-j}.
\end{eqnarray*}
Based on numerical verification, Lacasse presented the following conjecture. 
\begin{conjecture}
For any integer $n\geq 1$, there holds
\begin{eqnarray}\label{eqn 1.1}
\xi_2(n)\hskip-.22cm &=&\hskip-.22cm \xi(n)+n.
\end{eqnarray}
\end{conjecture}
Recently, by applying the Hurwitz identity on multivariate Abel polynomials, Younsi \cite{Younsi}
gave an algebraic proof of this conjecture. Later, using a decomposition of triply rooted trees into three
doubly rooted trees, Chen, Peng and Yang \cite{Chen} gave it a nice combinatorial interpretation.

In this note, using the derangement polynomials and their umbral representation, we provide another simple proof of (\ref{eqn 1.1}).

{\bf \section{The derangement polynomials and the proof of (\ref{eqn 1.1})} }

Recall that the derangement polynomials $\{\mathcal{D}_{n}(\lambda)\}_{n\geq 0}$ are
defined by
\begin{eqnarray}\label{eqn 2.1}
\mathcal{D}_{n}(\lambda) \hskip-.22cm &=&\hskip-.22cm \sum_{k=0}^{n}\binom{n}{k}D_k\lambda^{n-k}.
\end{eqnarray}
where $\mathcal{D}_{n}(1)=n!$ and $\mathcal{D}_n(0)=D_n$ is the $n$-th derangement number, counting
permutations on $[n]=\{1,2,\dots, n\}$ with no fixed points. The
derangement polynomials $\mathcal{D}_{n}(\lambda)$, also called
$\lambda$-factorials of $n$, have been considerably investigated by
Eriksen, Freij and W$\ddot{a}$stlund \cite{Eriksen}, Sun and Zhuang
\cite{SunZhuang}. The derangement polynomials $\mathcal{D}_{n}(\lambda)$ have the following basic property \cite{Eriksen}
and an Abel-type formula \cite{SunZhuang},
\begin{eqnarray}
\mathcal{D}_{n}(\lambda+\mu) \hskip-.22cm &=&\hskip-.22cm \sum_{k=0}^{n}\binom{n}{k}\mathcal{D}_{k}(\lambda)\mu^{n-k}, \label{eqn 2.2}\\
\mathcal{D}_{n}(\lambda+\mu) \hskip-.22cm &=&\hskip-.22cm \sum_{k=0}^{n}\binom{n}{k}(\lambda+k)^{k}(\mu-k-1)^{n-k}, \label{eqn 2.3}
\end{eqnarray}
and obey the recursive relation \cite{SunZhuang},
\begin{eqnarray}\label{eqn 2.4}
\sum_{k=0}^{n}\binom{n}{k}\mathcal{D}_{k}(\lambda)\mathcal{D}_{n-k}(\mu+1) \hskip-.22cm &=&\hskip-.22cm (\lambda+\mu-1)^{n+1}+(n-\lambda-\mu+2)\mathcal{D}_n(\lambda+\mu).
\end{eqnarray}
Denote by $\mathbf{D}$ the umbral operator defined by $\mathbf{D}^n=D_n$ for $n\geq 0$
(See \cite{Gesselb, Roman, RomRota} for more information on the umbral calculus), then by (\ref{eqn 2.1}) $\mathcal{D}_{n}(\lambda)$ can be represented as
\begin{eqnarray*}
\mathcal{D}_{n}(\lambda)  \hskip-.22cm &=&\hskip-.22cm (\mathbf{D}+\lambda)^n.
\end{eqnarray*}

Setting $\lambda=0, \mu=n+1$ in (\ref{eqn 2.4}), we have
\begin{eqnarray*}
n^{n+1}+\mathcal{D}_{n}(n+1)
\hskip-.22cm &=&\hskip-.22cm  \sum_{k=0}^{n}\binom{n}{k}\mathcal{D}_{k}(0)\mathcal{D}_{n-k}(n+2) \\
\hskip-.22cm &=&\hskip-.22cm  \sum_{k=0}^{n}\binom{n}{k}\mathcal{D}_{k}(0)\mu^{n-k}|_{\mu=(\mathbf{D}+n+2)} \\
\hskip-.22cm &=&\hskip-.22cm  \mathcal{D}_{n}(\mu)|_{\mu=(\mathbf{D}+n+2)}  \hskip2cm \mbox{by (\ref{eqn 2.2})}  \\
\hskip-.22cm &=&\hskip-.22cm  \sum_{k=0}^{n}\binom{n}{k}k^{k}(\mu-k-1)^{n-k}|_{\mu=(\mathbf{D}+n+2)}  \hskip1cm \mbox{by (\ref{eqn 2.3})}  \\
\hskip-.22cm &=&\hskip-.22cm  \sum_{k=0}^{n}\binom{n}{k}k^k\mathcal{D}_{n-k}(n-k+1),
\end{eqnarray*}
which proves (\ref{eqn 1.1}), if one notices that $\xi(n)$ and $\xi_2(n)$, by (\ref{eqn 2.3}), can be rewritten as
\begin{eqnarray*}
\xi(n)=\frac{1}{n^n}\mathcal{D}_{n}(n+1) \hskip.2cm \mbox{and}\hskip.2cm
\xi_2(n)= \frac{1}{n^n}\sum_{k=0}^{n}\binom{n}{k}k^k\mathcal{D}_{n-k}(n-k+1).
\end{eqnarray*}

\begin{remark}
By the nontrivial property of $\mathbf{D}$ \cite{SunZhuang},
\begin{eqnarray*}
(\mathbf{D}+\lambda)(\mathbf{D}+\lambda+n+1)^{n}\hskip-.22cm &=&\hskip-.22cm (n+\lambda)^{n+1},
\end{eqnarray*}
one can get another expression for $\xi_2(n)$,
\begin{eqnarray*}
\xi_2(n) \hskip-.22cm &=&\hskip-.22cm  \xi(n)+n= \frac{1}{n^n}\big(\mathcal{D}_{n}(n+1)+n^{n+1}\big) \\
         \hskip-.22cm &=&\hskip-.22cm  \frac{1}{n^n}\big((\mathbf{D}+n+1)^n+\mathbf{D}(\mathbf{D}+n+1)^n \big) \\
         \hskip-.22cm &=&\hskip-.22cm  \frac{1}{n^n}\big((\mathbf{D}+1)(\mathbf{D}+n+1)^n \big) =\frac{1}{n^n}\sum_{k=0}^{n}\binom{n}{k}(\mathbf{D}+1)^{k+1}n^{n-k}  \\
         \hskip-.22cm &=&\hskip-.22cm  \frac{1}{n^n}\sum_{k=0}^{n}\binom{n}{k}\mathcal{D}_{k+1}(1)n^{n-k}=\frac{1}{n^n}\sum_{k=0}^{n}\binom{n}{k}(k+1)!n^{n-k}.
\end{eqnarray*}
This expression has been obtained by Younsi using an identity of Hurwitz on multivariate Abel polynomials and plays a critical role in his proof.
\end{remark}

\vskip1cm

\section*{Acknowledgements} {  The work was partially supported by The National Science Foundation of China
and by the Fundamental Research Funds for the Central Universities.}

%==============================================================================================================

\end{document}